\documentclass[12pt]{article}
\usepackage{mathrsfs}
\usepackage{amssymb} \textwidth 150mm \textheight 220mm \oddsidemargin
15pt \evensidemargin 0pt \topmargin 0cm \headsep 0.3cm

\usepackage{amsmath}
\usepackage{graphicx}
\usepackage{indentfirst}
\usepackage{cite}
\usepackage{float}
\usepackage{CJK}

\begin{document}
\title{\Large\bf{On the number of zeros of Abelian integrals for discontinuous quadratic differential systems
\thanks{Corresponding author. E-mail addresses: zhaoliqin@bnu.edu.cn (L. Zhao), yangjh@mail.bnu.edu.cn (J. Yang), suisy@mail.bnu.edu.cn} }}
\author{{Jihua Yang$^{a,b}$, Liqin Zhao$^{a,*}$, Shiyou Sui$^a$ }
\\ {\small \it$^a$School of Mathematical Sciences, Beijing Normal University, Beijing 100875, China}\\
 {\small \it $^b$School of Mathematics and Computer Science, Ningxia Normal University, Guyuan}\\
 {\small \it   756000, China}}
\date{}
\maketitle \baselineskip=0.9\normalbaselineskip \vspace{-3pt}
\noindent
{\bf Abstract}\, Applying the Picard-Fuchs equation to the discontinuous differential system, we obtain the upper bounds of the number of zeros for Abelian integrals of four kinds of quadratic differential systems when it is perturbed inside all discontinuous polynomials with degree $n$. Furthermore, by using the {\it Chebyshev criterion}, we obtain the sharp upper bounds on each period annulus for $n=2$.
\vskip 0.2 true cm
\noindent
 \section{Introduction and main results}
 \setcounter{equation}{0}
\renewcommand\theequation{1.\arabic{equation}}

One of the main problems in the qualitative theory of continuous differential systems is the study of their limit cycles. The limit cycle of continuous differential systems has been studied intensively and many methodologies have been developed, such as Picard-Fuchs equation mentod \cite{GI,HI,YZ,ZZ}, Melnikov function method \cite{ANZ,GLT,TH}, averaging method \cite{LZL,LRR}, Chebyshev criterion \cite{GLLZ,GLT12,GMV,MV}, argument principle \cite{BL,LZL,P87,P89}.

It is common knowledge that the quadratic field with center must be one of the four types ($z=x+iy$)\cite{Z,I}:
 \vskip 0.2 true cm
$\dot{z}=-iz-z^2+2|z|^2+(b+ic)\bar{z}^2$,\ \, \, Hamiltonian $Q_3^H$.\vskip 0.2 true cm

 $\dot{z}=-iz+az^2+2|z|^2+b\bar{z}^2$,\quad\qquad\ Reversible $Q_3^R$.\vskip 0.2 true cm

 $\dot{z}=-iz+4z^2+2|z|^2+(b+ic)\bar{z}^2$,\ \ $|b+ic|=2$, codimension four $Q_4$. \vskip 0.2 true cm

 $\dot{z}=-iz+z^2+(b+ic)\bar{z}^2$,\quad\qquad\ \ \ \ generalized Lotka-Volterra $Q_3^{LV}$.\vskip 0.2 true cm
 The number of limit cycles of these system has been investigated in hundreds of excellent papers, see \cite{HXC,HXM,GGI,LZLZ,SZ,XHX,ZLLZ} and the references therein.  The reversible type $Q_3^R$ is equivalent to
$$\dot{x}=y+(a+b+2)x^2-(a+b-2)y^2,\ \ \dot{y}=-x+2(a-b)xy.$$
Taking the linear transformation $x_1=1-2(a-b)y$, $y_1=x$, $dt_1=-2(a-b)dt$, we obtain the new system
 \begin{eqnarray}
 \dot{x}=-xy,\ \ \dot{y}=-\frac{a+b+2}{2(a-b)}y^2+\frac{a+b-2}{8(a-b)^3}x^2
 -\frac{b-1}{2(a-b)^3}x-\frac{a-3b+2}{8(a-b)^3},
 \end{eqnarray}
 which has the first integral
 \begin{eqnarray}
 H(x,y)=x^{-\frac{a+b+2}{a-b}}\Big[\frac{1}{2}y^2+
 \frac{1}{8(a-b)^2}\Big(\frac{a+b-2}{a-3b-2}x^2+2\frac{b-1}{b+1}x+
 \frac{a-3b+2}{a+b+2}\Big)\Big],
 \end{eqnarray}
 and the integrating factor $\mu(x,y)=-x^{-2\frac{a+1}{a-b}}$ (here and below, we shall omit the subscript 1).
Taking $(a,b)=(-11,-3)$ and $(5,-3)$, we have (r19) and (r20) as follows:
\begin{eqnarray}
\dot{x}=-xy,\ \ \dot{y}=-\frac{3}{4}y^2+\frac{1}{2^8}x^2-\frac{1}{2^8}x,
\end{eqnarray}
\begin{eqnarray}
\dot{x}=-xy,\ \ \dot{y}=-\frac{1}{4}y^2+\frac{1}{2^8}x-\frac{1}{2^8},
\end{eqnarray}
whose first integrals are
\begin{eqnarray}
H(x,y)=x^{-\frac{3}{2}}\Big(\frac{1}{2}y^2+\frac{1}{2^7}x^2+\frac{1}{2^7}x\Big)
=h,\ \ h\in\Big(\frac{1}{2^6},+\infty\Big)
\end{eqnarray}
with the integrating factor $\mu(x,y)=-x^{-\frac{5}{2}}$ and
\begin{eqnarray}
H(x,y)=x^{-\frac{1}{2}}\Big(\frac{1}{2}y^2+\frac{1}{2^7}x+\frac{1}{2^7}\Big)
=h,\ \ h\in\Big(\frac{1}{2^6},+\infty\Big)
\end{eqnarray}
with the integrating factor $\mu(x,y)=-x^{-\frac{3}{2}}$, respectively.

For a quadratic differential system of the form
\begin{eqnarray}
\dot{x}=-y+X_2(x,y),\quad \dot{y}=x+Y_2(x,y),
\end{eqnarray}
where $X_2$, $Y_2$ are quadratic homogeneous polynomials, the origin is a center, Loud \cite{L} first classified system (1.7) with an isochronous center: The origin is an isochronous center of system (1.7) if and only if the system can be brought to one of the following systems $S_1$, $S_2$, $S_3$ and $S_4$, through a linear change of coordinates and a rescaling of time, where
\begin{eqnarray*}
\begin{aligned}
&S_1:\ \dot{x}=-y+x^2-y^2,\ \ \ \dot{y}=x+2xy,\\
&S_2:\ \dot{x}=-y+x^2,\ \ \ \dot{y}=x+xy,\\
&S_3:\ \dot{x}=-y-\frac{4}{3}x^2,\ \ \ \dot{y}=x-\frac{16}{3}xy,\\
&S_4:\ \dot{x}=-y+\frac{16}{3}x^2-\frac{4}{3}y^2,\ \ \ \dot{y}=x+\frac{8}{3}xy.\\
\end{aligned}
\end{eqnarray*}
Taking linear transformations $x_2=1+2y$, $y_2=\sqrt{2}x$ and $x_2=1+y$, $y_2=\sqrt{2}x$ for $S_1$ and $S_2$ respectively, we obtain
\begin{eqnarray}
        \dot{x}=\sqrt{2}xy, \ \ \
        \dot{y}=\frac{\sqrt{2}}{4}(1-x^2+2y^2)
       \end{eqnarray}
and
\begin{eqnarray}
        \dot{x}=\frac{\sqrt{2}}{2}xy, \ \ \
        \dot{y}=\frac{\sqrt{2}}{2}(2-2x+y^2)
        \end{eqnarray}
(here and below, we shall omit the subscript 2), whose first integrals are
\begin{eqnarray}
H(x,y)=x^{-1}\Big(\frac{1}{2}y^2+\frac{1}{4}x^2-\frac{1}{2}x+\frac{1}{4}\Big)=h,\ h\in(-\infty,-1)\cup(0,+\infty)
\end{eqnarray}
with integrating factor $\mu(x,y)=\frac{\sqrt{2}}{2}x^{-2}$
and
\begin{eqnarray}
H(x,y)=x^{-2}\Big(\frac{1}{2}y^2+x^2-2x+1\Big)=h,\ h\in(0,1)
\end{eqnarray}
with integrating factor $\mu(x,y)=\sqrt{2}x^{-3}$, respectively.

Chicone and Jacobs proved that in \cite{CJ}, under all continuous quadratic polynomial perturbations, at most 2 limit cycles bifurcates from the period annuli of $S_1$, $S_2$, $S_3$ and $S_4$. Their study is based in the displacement function using some results of Bautin \cite{B}. Li et al.\cite{LLL} studied the quadratic isochronous centers $S_1$, $S_2$, $S_3$ and $S_4$ under perturbations inside the class of all continuous polynomial systems of degree $n$. They obtained the exact upper bounds of the number of limit cycles for $S_1$, $S_2$ and $S_3$ and an upper bound for $S_4$ by analyzing the corresponding Melnikov function of first order.

Stimulated by discontinuous phenomena in the real world, such as biology \cite{K}, nonlinear oscillations \cite{T}, impact and friction mechanics \cite{BBC}, a big interest has appeared for studying the number of limit cycles and their relative positions of discontinuous differential systems. This problem can be seen as an extension of the infinitesimal Hilbert's 16th Problem  to the discontinuous world.

As far as we know, in addition to the Picard-Fuchs equation method, the Melnikov function method, averaging method, Chebyshev criterion and argument principle have been used to study the bifurcations of limit cycle of discontinuous differential systems, see \cite{YZ16,LHR,LH,YZ17,LMN,CLZ,LC} and the references therein. By the averaging method of first order, Llibre and Mereu \cite{LM} studied the number of limit cycles bifurcated from the period annuli of quadratic isochronous system $S_1$ and $S_2$ under perturbations of discontinuous quadratic polynomials, they found that at least 4 (resp. 5) limit cycles can bifurcate from the isochronous center at the origin of $S_1$ (resp. $S_2$). By using the averaging method and Chebyshev criterion, Li and Cen \cite{LC} studied the quadratic isochronous system $S_3$, and Cen, Li and Zhao \cite{CLZ} studied the quadratic isochronous system $S_4$. They found that at most 4 limit cycles can bifurcated from the isochronous center at the origin of $S_3$, and at most 5 limit cycles can bifurcated from the isochronous center at the origin of $S_4$ and both of them can be realizable under perturbations of discontinuous quadratic polynomials.

In this paper, by using the Picard-Fuchs method to discontinuous differential systems, we investigate the number of zeros of Abelian integrals for the systems $(r19)$, $(r20)$, $S_1$ and $S_2$ under arbitrary discontinuous polynomial perturbations of degree $n$. For the sake of convenience, we rewrite the first integrals (1.5), (1.6), (1.10) and (1.11) as follows:
\begin{eqnarray}
H(x,y)=x^{-k}\Big(\frac{1}{2}y^2+\lambda_2x^2+\lambda_1x+\lambda_0\Big), \ h\in\Sigma,
\end{eqnarray}
where

$\bullet$ $S_1$: $k=1$, $\lambda_2=\frac{1}{4}$, $\lambda_1=-\frac{1}{2}$, $\lambda_0=\frac{1}{4}$, $\Sigma=(-\infty,-1)\cup(0,+\infty)$.
 \vskip 0.2 true cm
$\bullet$ $S_2$: $k=2$, $\lambda_2=1$, $\lambda_1=-2$, $\lambda_0=1$, $\Sigma=(0,1)$.
 \vskip 0.2 true cm
$\bullet$ $(r19)$: $k=\frac{3}{2}$, $\lambda_2=\frac{1}{2^7}$, $\lambda_1=\frac{1}{2^7}$, $\lambda_0=0$, $\Sigma=(\frac{1}{2^6},+\infty)$.
 \vskip 0.2 true cm
 $\bullet$ $(r20)$: $k=\frac{1}{2}$, $\lambda_2=0$, $\lambda_1=\frac{1}{2^7}$, $\lambda_0=\frac{1}{2^7}$, $\Sigma=(\frac{1}{2^6},+\infty)$.
 \vskip 0.2 true cm

 The system (1.8) has two isochronous centers ($\pm$1,0) and $h=0$ and $h=-1$ correspond to (1,0) and (-1,0) respectively, the system (1.8) has an isochronous center (1,0) and $h=0$ corresponds to (1,0), the system (1.3) has an elementary center $(1,0)$ and $h=\frac{1}{2^6}$ corresponds to (1,0), and the system (1.4) has an elementary center (1,0) and $h=\frac{1}{2^6}$ corresponds to (1,0), see Figs.\,1-4.
\begin{figure}[htbp]
 \centering \includegraphics[width=3in]{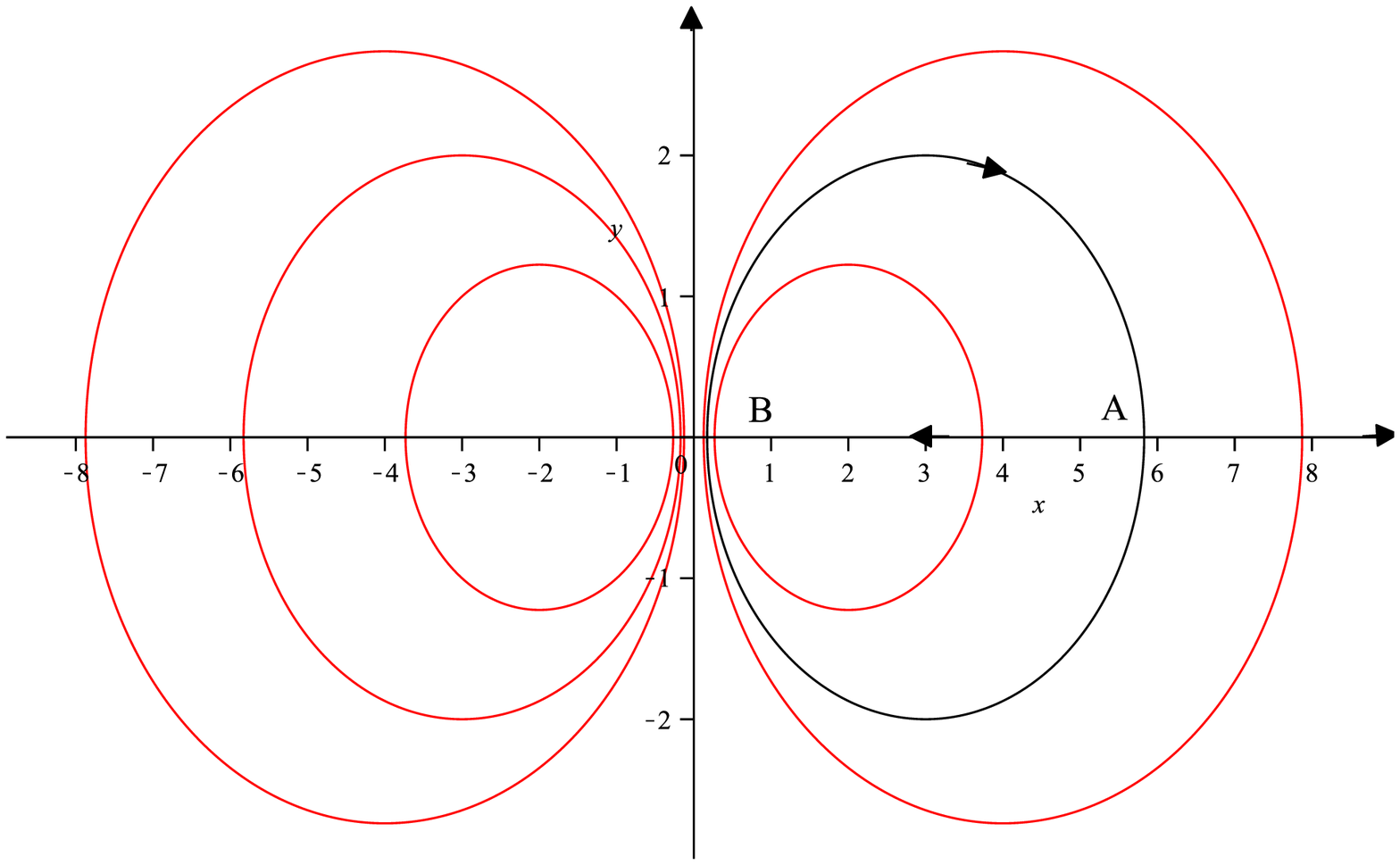}
\begin{center}
{\small{\bf Fig.\,1.} The phase portrait of system (1.8). }
\end{center}
\end{figure}
\begin{figure}[htbp]
 \centering \includegraphics[width=3in]{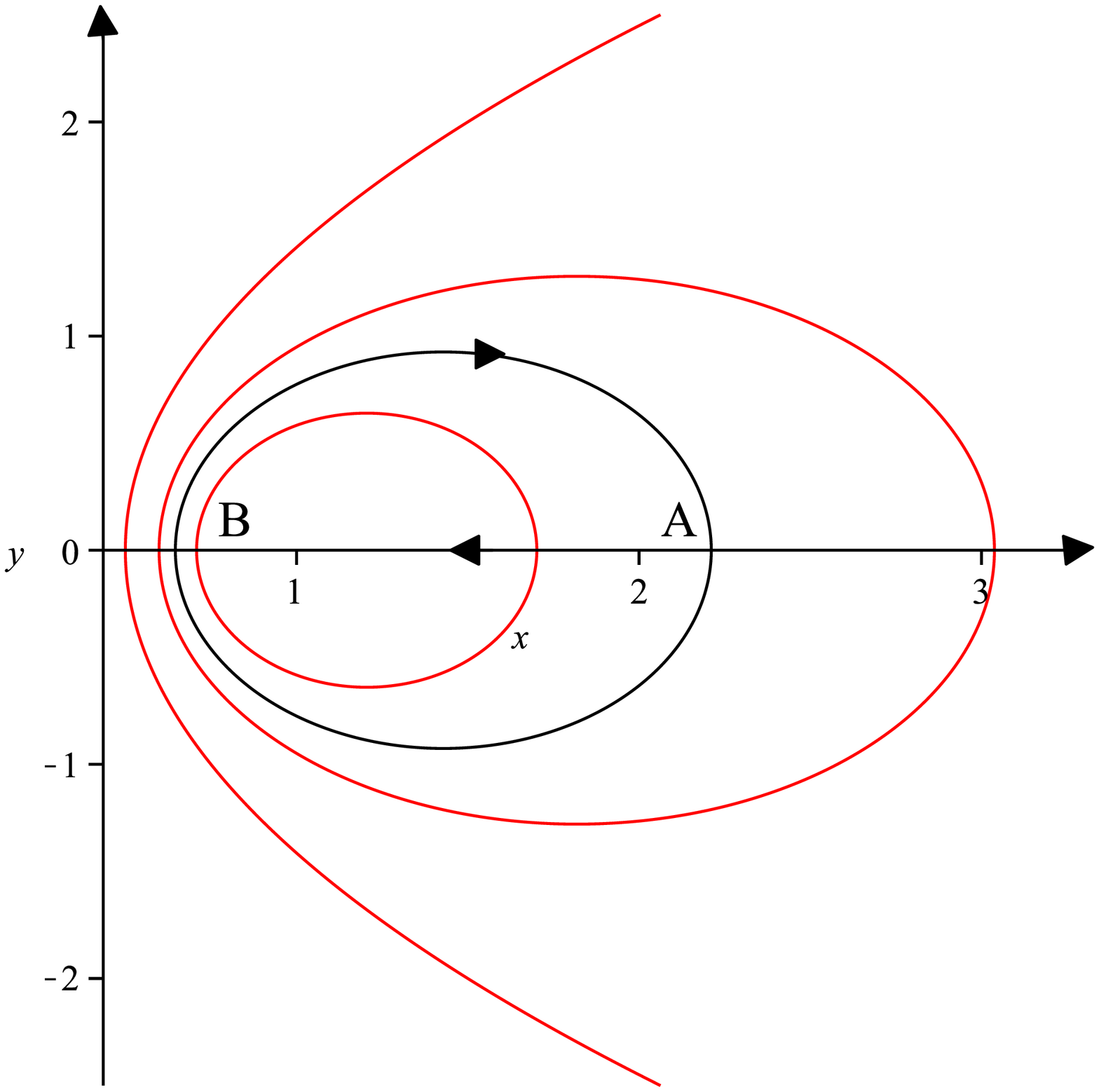}
\begin{center}
{\small{\bf Fig.\,2.} The phase portrait of system (1.9). }
\end{center}
\end{figure}
\begin{figure}[htbp]
 \centering \includegraphics[width=3in]{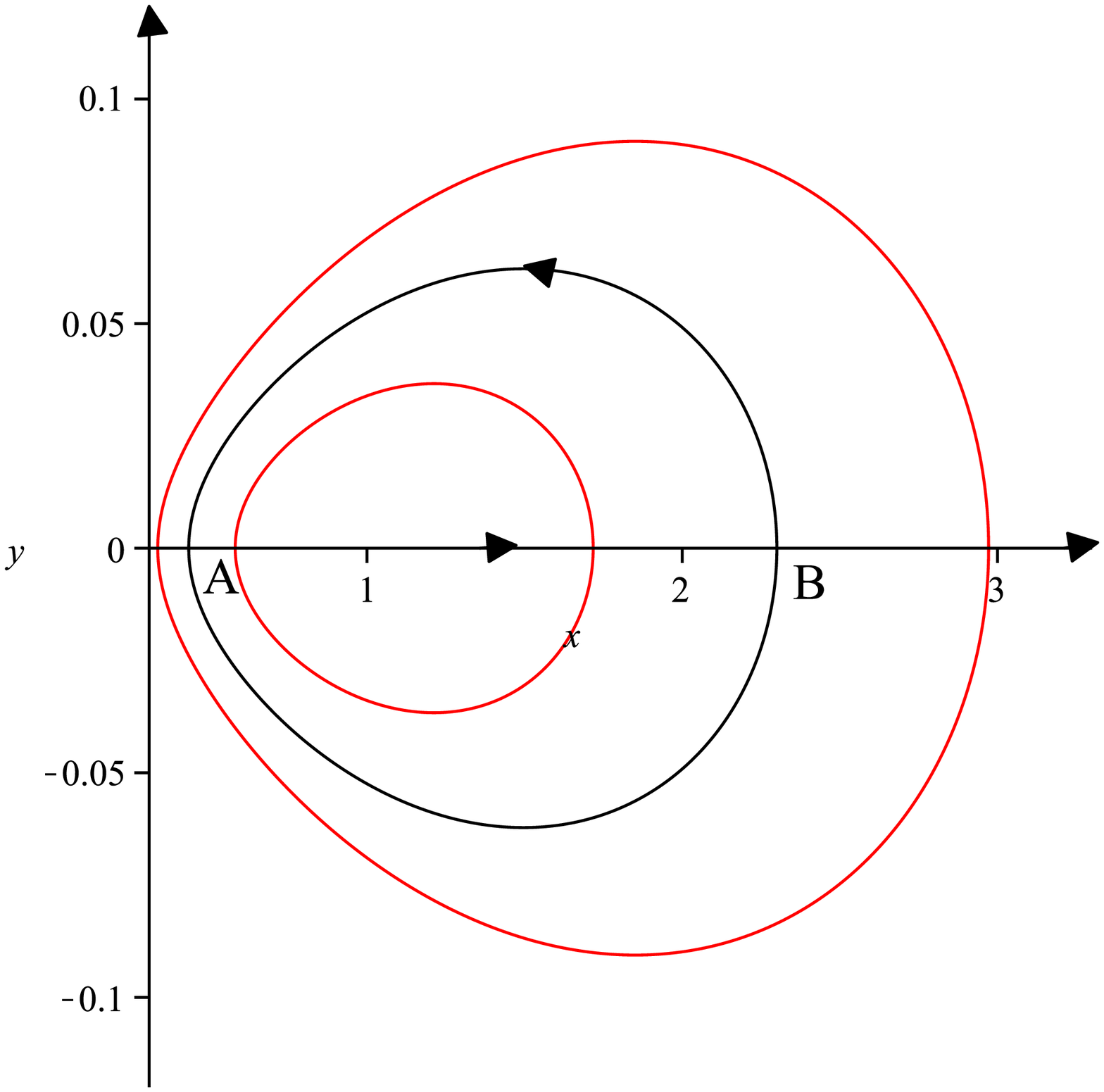}
\begin{center}
{\small{\bf Fig.\,3.} The phase portrait of system (1.3). }
\end{center}
\end{figure}
\begin{figure}[htbp]
 \centering \includegraphics[width=3in]{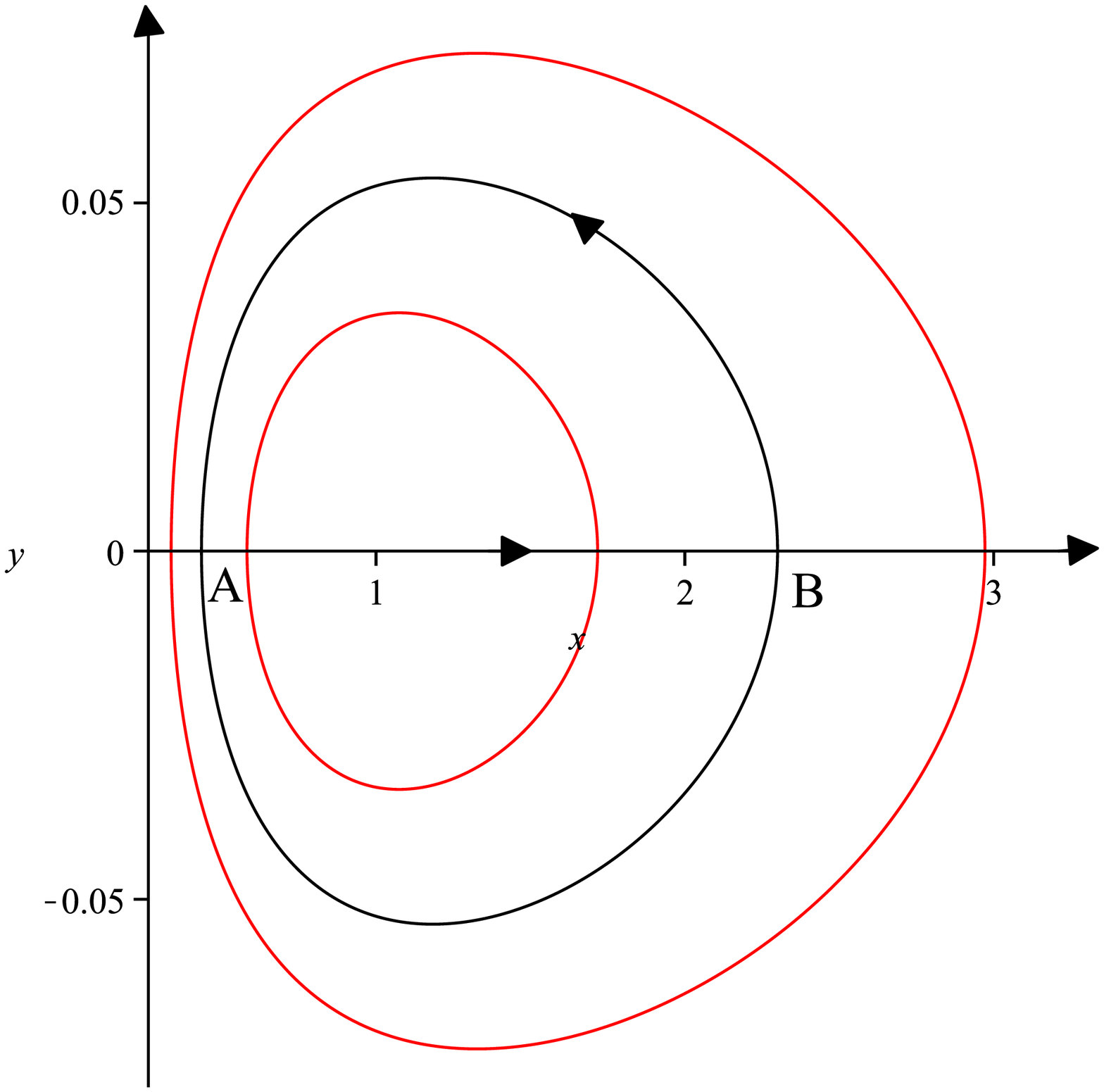}
\begin{center}
{\small{\bf Fig.\,4.} The phase portrait of system (1.4). }
\end{center}
\end{figure}

The perturbed systems of (1.8), (1.9), (1.3) and (1.4) are
\begin{eqnarray}
\left(
  \begin{array}{c}
          \dot{x} \\
          \dot{y}
 \end{array}
 \right)=\begin{cases}
 \left(
  \begin{array}{c}
          \sqrt{2}xy+\varepsilon f^+(x,y) \\
          \frac{\sqrt{2}}{4}(1-x^2+2y^2)+\varepsilon g^+(x,y)
 \end{array}
 \right), \quad y>0,\\
 \left(
  \begin{array}{c}
         \sqrt{2}xy+\varepsilon f^-(x,y) \\
          \frac{\sqrt{2}}{4}(1-x^2+2y^2)+\varepsilon g^-(x,y)
 \end{array}
 \right),\quad y<0,\\
 \end{cases}
       \end{eqnarray}
\begin{eqnarray}
\left(
  \begin{array}{c}
          \dot{x} \\
          \dot{y}
 \end{array}
 \right)=\begin{cases}
 \left(
  \begin{array}{c}
          \frac{\sqrt{2}}{2}xy+\varepsilon f^+(x,y) \\
          \frac{\sqrt{2}}{2}(2-2x+y^2)+\varepsilon g^+(x,y)
 \end{array}
 \right), \quad y>0,\\
 \left(
  \begin{array}{c}
         \frac{\sqrt{2}}{2}xy+\varepsilon f^+(x,y) \\
          \frac{\sqrt{2}}{2}(2-2x+y^2)+\varepsilon g^+(x,y)
 \end{array}
 \right),\quad y<0,\\
 \end{cases}
       \end{eqnarray}
 \begin{eqnarray}
\left(
  \begin{array}{c}
          \dot{x} \\
          \dot{y}
 \end{array}
 \right)=\begin{cases}
 \left(
  \begin{array}{c}
          -xy+\varepsilon f^+(x,y) \\
          -\frac{3}{4}y^2+\frac{1}{2^8}x^2-\frac{1}{2^8}x+\varepsilon g^+(x,y)
 \end{array}
 \right), \quad y>0,\\
 \left(
  \begin{array}{c}
         -xy+\varepsilon f^-(x,y) \\
          -\frac{3}{4}y^2+\frac{1}{2^8}x^2-\frac{1}{2^8}x+\varepsilon g^-(x,y)
 \end{array}
 \right),\quad y<0\\
 \end{cases}
       \end{eqnarray}
       and
\begin{eqnarray}
\left(
  \begin{array}{c}
          \dot{x} \\
          \dot{y}
 \end{array}
 \right)=\begin{cases}
 \left(
  \begin{array}{c}
          -xy+\varepsilon f^+(x,y) \\
          -\frac{1}{4}y^2+\frac{1}{2^8}x-\frac{1}{2^8}+\varepsilon g^+(x,y)
 \end{array}
 \right), \quad y>0,\\
 \left(
  \begin{array}{c}
         -xy+\varepsilon f^+(x,y) \\
          -\frac{1}{4}y^2+\frac{1}{2^8}x-\frac{1}{2^8}+\varepsilon g^+(x,y)
 \end{array}
 \right),\quad y<0,\\
 \end{cases}
       \end{eqnarray}
respectively, where $0<|\varepsilon|\ll1$,
$$f^\pm(x,y)=\sum\limits_{i+j=0}^na^\pm_{i,j}x^iy^j,\ \ g^\pm(x,y)=\sum\limits_{i+j=0}^nb^\pm_{i,j}x^iy^j.$$

From the Proposition 1.1 in \cite{YZ17}, we know that the number of zeros of the following Abelian integral
\begin{eqnarray}
\begin{aligned}
M(h)=&\int_{\Gamma^+_h}x^{-k-1}[g^+(x,y)dx-f^+(x,y)dy]\\
&+\int_{\Gamma^-_h}x^{-k-1}[g^-(x,y)dx-f^-(x,y)dy],\ \  h\in\Sigma,
\end{aligned}
\end{eqnarray}
controls the number of limit cycles of the perturbed systems (1.13)-(1.16) for $|\varepsilon|$ small enough,
where
$$\begin{aligned}&\Gamma_h^+=\{(x,y)|H(x,y)=h,h\in\Sigma,y>0\},\\
&\Gamma^-_h=\{(x,y)|H(x,y)=h,h\in\Sigma,y<0\}.\end{aligned}$$

Our main results are the following two theorems.
 \vskip 0.2 true cm

\noindent
{\bf Theorem 1.1.}\, {\it Suppose that $h\in\Sigma$.}
 \vskip 0.2 true cm
\noindent
 (i) {\it For $S_1$, the upper bound of the number of zeros of the Abelian integrals is  $4n+30$ if $n\geq0$ counting multiplicity. For $S_2$, the upper bound of the number of zeros of the Abelian integrals is  $10n-4$ if $n\geq1$, 2 if $n=0$ counting multiplicity.}
 \vskip 0.2 true cm
 \noindent
 (ii) {\it For $(r19)$, the upper bound of the number of zeros of the Abelian integrals is  $4n-3$ if $n\geq4$, 11 if $n=0,1,2,3$ counting multiplicity. For $(r20)$, the upper bound of the number of zeros of the Abelian integrals is  $4n+3$ if $n\geq3$, 8 if $n=0,1,2$ counting multiplicity.}

 \vskip 0.2 true cm

\noindent
{\bf Theorem 1.2.}\, {\it Suppose that $h\in\Sigma$, $a^+_{i,j}=a^-_{i,j}$ and $b^+_{i,j}=b^-_{i,j}$.}
 \vskip 0.2 true cm
\noindent
 (i) {\it For $S_1$, the upper bound of the number of zeros of the Abelian integrals is  $2n$ if $n\geq0$ counting multiplicity. For $S_2$, the upper bound of the number of zeros of the Abelian integrals is  $2n-1$ if $n\geq1$, 1 if $n=0$ counting multiplicity. }
 \vskip 0.2 true cm
 \noindent
 (ii) {\it For $(r19)$, the upper bound of the number of zeros of the Abelian integrals is  $2n-3$ if $n\geq4$, 4 if $n=0,1,2,3$ counting multiplicity. For $(r20)$, the upper bound of the number of zeros of the Abelian integrals is  $2n$ if $n\geq3$, 3 if $n=0,1,2$ counting multiplicity.}

\vskip 0.2 true cm

\noindent
{\bf Theorem 1.3.}\, {\it For $n=2$, the exact number of zeros of the Abelian integrals associated to $S_1$ (resp. $S_2$) is 5 (resp. 6) counting multiplicity for $h\in(-\infty,-1)$ or $h\in(0,+\infty)$ (resp. $h\in(0,1)$). }

 \vskip 0.2 true cm

\noindent
{\bf Remark 1.1.}\, (i) It is easy to check that $S_1\in Q_3^{LV}$, $S_2\in Q_3^R$.
 \vskip 0.2 true cm

\noindent
(ii) If $a^+_{i,j}=a^-_{i,j}$ and $b^+_{i,j}=b^-_{i,j}$, X. Hong et al.\cite{HXM,HXC} obtain that the upper bound of the number of zeros of the Abelian integrals is  $6n-12$ if $n\geq4$, 10 if $n=0,1,2,3$ for $(r19)$ and $6n-9$ if $n\geq3$, 7 if $n=0,1,2$ for $(r20)$, counting multiplicity.

 \vskip 0.2 true cm

\noindent
(ii) In \cite{LM}, Llibre and Mereu proved that $S_1$ and $S_2$ under the perturbations of discontinuous quadratic polynomials without constant terms have at least 4 and 5 limit cycles respectively.

\section{The algebraic structure of $M(h)$ and Picard-Fuchs equation}
 \setcounter{equation}{0}
\renewcommand\theequation{2.\arabic{equation}}

In this section, we obtain the algebraic structure of Abelian integral $M(h)$.
For $h\in\Sigma$, we denote
\begin{eqnarray}
I_{i,j}(h)=\int_{\Gamma^+_h}x^{i-k-1}y^jdx,\ \ J_{i,j}(h)=\int_{\Gamma^-_h}x^{i-k-1}y^jdx,
\end{eqnarray}
where $i$ and $j$ are integers. We first prove the following results.

 \vskip 0.2 true cm

\noindent
{\bf Lemma 2.1.}\, {\it Suppose that $h\in\Sigma$.}
\vskip 0.2 true cm

\noindent
(i) {\it For $S_1$,
\begin{eqnarray}
\begin{aligned}
M(h)=\frac{1}{2h+1}[\alpha_1(h) I_{0,0}(h)+\beta_1(h) I_{1,1}(h)+\gamma_1(h) I_{-1,1}(h)+\delta_1(h) I_{0,2}(h)].
\end{aligned}
\end{eqnarray}
where $\alpha_1(h)$, $\beta_1(h)$, $\gamma_1(h)$ and $\delta_1(h)$ are polynomials of $h$ with
\begin{eqnarray*}\begin{aligned}
\deg \alpha_1(h), \deg\beta_1(h)\leq n-1,\ \deg\gamma_1(h)\leq n-3,\ \delta_1(h)\leq2.
\end{aligned}\end{eqnarray*}}

\noindent
(ii) {\it  For $S_2$,
\begin{eqnarray}
\begin{aligned}
M(h)=\frac{1}{(h-1)^{n-2}}[\alpha_2(h) I_{0,1}(h)+\beta_2(h) I_{1,0}(h)+\gamma_2(h) I_{1,1}(h)+\delta_2(h) I_{0,2}(h)].
\end{aligned}
\end{eqnarray}
 where $\alpha_2(h)$, $\beta_2(h)$, $\gamma_2(h)$ and $\delta_2(h)$ are polynomials of $h$ with
\begin{eqnarray*}\begin{aligned}
\deg \alpha_2(h), \deg\delta_2(h)\leq n-2,\ \deg\beta_2(h),\deg\gamma_2(h)\leq n-1.
\end{aligned}\end{eqnarray*}}

\noindent
(iii) {\it For $(r19)$,
\begin{eqnarray}
M(h)=\tau_3h\ln (64h-\sqrt{4096h^2-1})+\alpha_3(h)I_{\frac{1}{2},1}
+\beta_3(h)I_{1,1}+\gamma_3(h)I_{1,0},
\end{eqnarray}
where $\tau_3$ is a constant, and $\alpha_3(h)$, $\beta_3(h)$ and $\gamma_3(h)$ are polynomials of $h$ with
$$\begin{aligned}
&\deg\alpha_3(h)\leq2n-5,\ \deg\beta_3(h),\deg\gamma_3(h)\leq2n-4,\ n\geq4,\\
&\deg\alpha_3(h)\leq3,\ \deg\beta_3(h),\deg\gamma_3(h)\leq2,\ n=0,1,2,3.\end{aligned}$$}

\noindent
(iv) {\it For $(r20)$,
\begin{eqnarray}
M(h)=\tau_4h\ln (64h-\sqrt{4096h^2-1})+\alpha_4(h)I_{0,1}+\beta_4(h)I_{\frac{1}{2},1}+\gamma_4(h)I_{0,0},
\end{eqnarray}
where $\tau_4$ is a constant, and $\alpha_4(h)$, $\beta_4(h)$ and $\gamma_4(h)$ are polynomials of $h$ with
\begin{eqnarray*}\begin{aligned}&\deg\alpha_4(h)\leq2n-4,\ \deg\beta_4(h)\leq2n-3,\ \deg\gamma_4(h)\leq2n-2,\ n\geq3,\\
&\deg\alpha_4(h)\leq2,\ \deg\beta_4(h)\leq1,\ \deg\gamma_4(h)\leq2,\ n=0,1,2.\end{aligned}\end{eqnarray*}}

 \vskip 0.2 true cm

\noindent
{\bf Proof.}\, Let $D$ be the interior of $\Gamma_{h}^+\cup \overrightarrow{AB}$, see the black line in Figs.\,1-4. Using the Green's Formula, we have for $j\geq0$
\begin{eqnarray*}
\begin{aligned}
\int_{\Gamma^+_{h}}x^iy^jdy
=&\oint_{\Gamma^+_{h}\cup\overrightarrow{AB}}x^iy^jdy-\int_{\overrightarrow{AB}}x^iy^jdy\\
=&\oint_{\Gamma^+_{h}\cup\overrightarrow{AB}}x^iy^jdy
=\mp i\iint\limits_Dx^{i-1}y^jdxdy,
\end{aligned}
\end{eqnarray*}
\begin{eqnarray*}
\begin{aligned}
\int_{\Gamma^+_{h}}x^{i-1}y^{j+1}dx=
\oint_{\Gamma^+_{h}\cup\overrightarrow{AB}}x^{i-1}y^{j+1}dx
=\pm(j+1)\iint\limits_Dx^{i-1}y^jdxdy,
\end{aligned}
\end{eqnarray*}
where $\pm$ is determined by the direction of the integral path $\Gamma^+_{h}\cup\overrightarrow{AB}$.
Hence, \begin{eqnarray}\int_{\Gamma^+_{h}}x^iy^jdy=-\frac{i}{j+1}\int_{\Gamma^+_{h}}x^{i-1}y^{j+1}dx,\ j\geq0.\end{eqnarray}
In a similar way, we have
\begin{eqnarray}\int_{\Gamma^-_{h}}x^iy^jdy=-\frac{i}{j+1}\int_{\Gamma^-_{h}}x^{i-1}y^{j+1}dx,\ j\geq0.\end{eqnarray}
By a straightforward calculation and noting that (2.6) and (2.6), we obtain
\begin{eqnarray}
\begin{aligned}
M(h)=&\int_{\Gamma^+_h}x^{-k-1}\big(g^+(x,y)dx-f^+(x,y)dy\big)\\
&+\int_{\Gamma^-_h}x^{-k-1}\big(g^-(x,y)dx-f^-(x,y)dy)\\
=&\int_{\Gamma_h^+}\sum\limits_{i+j=0}^nb^+_{i,j}x^{i-k-1}y^jdx-\int_{\Gamma_h^+}\sum\limits_{i+j=0}^na^+_{i,j}x^{i-k-1}y^jdy\\
&+\int_{\Gamma_h^-}\sum\limits_{i+j=0}^nb^-_{i,j}x^{i-k-1}y^jdx-\int_{\Gamma_h^-}\sum\limits_{i+j=0}^na^-_{i,j}x^{i-k-1}y^jdy\\
=&\sum\limits_{i+j=0}^nb^+_{i,j}\int_{\Gamma^+_h}x^{i-k-1}y^jdx+
\sum\limits_{i+j=0}^n\frac{i-k-1}{j+1}a^+_{i,j}\int_{\Gamma^+_h}x^{i-k-2}y^{j+1}dx\\
&+\sum\limits_{i+j=0}^nb^-_{i,j}\int_{\Gamma^-_h}x^{i-k-1}y^jdx+
\sum\limits_{i+j=0}^n\frac{i-k-1}{j+1}a^-_{i,j}\int_{\Gamma^-_h}x^{i-k-2}y^{j+1}dx\\
=&\sum\limits_{i+j=0,i\geq-1,j\geq0}^n\tilde{a}_{i,j}I_{i,j}(h)
+\sum\limits_{i+j=0,i\geq-1,j\geq0}^n\tilde{b}_{i,j}J_{i,j}(h)\\
:=&\sum\limits_{i+j=0,i\geq-1,j\geq0}^n\rho_{i,j}I_{i,j}(h),
\end{aligned}
\end{eqnarray}
where in the last equality we have used that $I_{i,j}(h)=(-1)^{j+1}J_{i,j}(h)$.

Differentiating (1.12) with respect to $x$, we obtain
\begin{eqnarray}
-\frac{k}{2}x^{-k-1}y^2+x^{-k}y\frac{\partial y}{\partial x}+\lambda_2(2-k)x^{1-k}+\lambda_1(1-k)x^{-k}-\lambda_0kx^{-k-1}=0.
\end{eqnarray}
Multiplying (2.9) by $x^iy^{j-2}dx$, integrating over $\Gamma^+_h$ and noting that (2.6), we have
\begin{eqnarray}
(2i+kj-2k)I_{i,j}=2j\big[\lambda_2(2-k)I_{i+2,j-2}+\lambda_1(1-k)I_{i+1,j-2}-\lambda_0kI_{i,j-2}\big].
\end{eqnarray}
Similarly, multiplying (1.12) by $x^{i-k-1}y^jdx$ and integrating over $\Gamma^+_h$ yields
\begin{eqnarray}
hI_{i,j}=\frac{1}{2}I_{i-k,j+2}+\lambda_2I_{i-k+2,j}+\lambda_1I_{i-k+1,j}+\lambda_0I_{i-k,j}.
\end{eqnarray}

(i). For $S_1$, elementary manipulations reduce Eqs. (2.10) and (2.11) to
\begin{eqnarray}
I_{i,j}=\frac{1}{i+j-1}\big[(2i+j-4)(2h+1)I_{i-1,j}-(i-3)I_{i-2,j}\big]
\end{eqnarray}
and
\begin{eqnarray}
I_{i,j}=\frac{j}{2(i+j-1)}\big[(2h+1)I_{i+1,j-2}-I_{i,j-2}\big].
\end{eqnarray}
From (2.10) we have
\begin{eqnarray}
I_{2,0}=I_{0,0},\ \ I_{-1,3}=\frac{3}{2}(I_{-1,1}-I_{1,1}).
\end{eqnarray}
Taking $(i,j)=(1,1)$ in (2.12) and $(i,j)=(0,2)$ in (2.13), respectively, we have
\begin{eqnarray}
I_{1,1}=-(2h+1)I_{0,1}+2I_{-1,1},\ I_{0,2}=(2h+1)I_{1,0}-I_{0,0}.
\end{eqnarray}
Taking $(i,j)=(-1,2)$ in (2.12), $(i,j)=(0,0)$ in (2.13) and eliminating $I_{-1,0}$, we obtain
\begin{eqnarray}
I_{-1,2}=(h+\frac{1}{2})I_{0,0}-\frac{1}{2}I_{1,0}.
\end{eqnarray}
From (2.12) and (2.13) we obtain
\begin{eqnarray}
\begin{aligned}
&I_{-1,4}=(2h+1)I_{0,2}-I_{-1,2},\ \ I_{0,3}=\frac{3}{4}(2h+1)I_{1,1}-\frac{3}{4}I_{0,1},\\
&I_{1,2}=I_{-1,2},\ \ I_{2,1}=\frac{1}{2}(2h+1)I_{1,1}+\frac{1}{2}I_{0,1},\ \ I_{3,0}=(2h+1)I_{2,0}.
\end{aligned}
\end{eqnarray}

Now we prove (2.2) by induction on $n$. In fact, (2.14)-(2.17) imply that (2.2) holds for $n=2,3$. Now assume that statement (2.2) holds for $i+j\leq k-1\, (k\geq3)$, then for $i+j=k$ taking $(i,j)=(-1,k+1),(0,k),(1,k-1),\cdots,(k-2,2)$ in (2.13) and $(i,j)=(k-1,1),(k,0)$ in (2.12), respectively, we have
\begin{eqnarray}
\left(\begin{matrix}
                I_{-1,k+1}\\
                 I_{0,k}\\
                 I_{1,k-1}\\
                 \vdots\\
                  I_{k-2,2}\\
                  I_{k-1,1}\\
                  I_{k,0}
                \end{matrix}\right)\ \
=\frac{1}{k-1}\left(\begin{matrix}
                \frac{k+1}{2}\big[(2h+1)I_{0,k-1}-I_{-1,k-1}\big]\\
               \frac{k}{2}\big[(2h+1)I_{1,k-2}-I_{0,k-2}\big]\\
                \frac{k-1}{2}\big[(2h+1)I_{2,k-3}-I_{1,k-3}\big]\\
                                  \vdots\\
                  (2h+1)I_{k-1,0}-I_{k-2,0}\\
                (2k-5)(2h+1)I_{k-2,1}-(k-4)I_{k-3,1}\\
                (2k-4)(2h+1)I_{k-1,0}-(k-3)I_{k-2,0}
                \end{matrix}\right).
\end{eqnarray}
By the induction hypothesis we obtain the statement (2.2). From (2.18) we have
\begin{eqnarray*}
\begin{aligned}
I_{i,j}(h)=&\alpha^{(k-2)}(h)I_{0,0}+\beta^{(k-2)}(h)I_{1,1}+\gamma^{(k-2)}(h)I_{-1,1}+\delta^{(k-2)}(h)I_{0,2}\\
&+h\big[\alpha^{(k-1)}(h)I_{0,0}+\beta^{(k-1)}(h)I_{1,1}+\gamma^{(k-1)}(h)I_{-1,1}+\delta^{(k-1)}(h)I_{0,2}\big]\\
:=&\alpha^{(k)}(h)I_{0,0}+\beta^{(k)}(h)I_{1,1}+\gamma^{(k)}(h)I_{-1,1}+\delta^{(k)}(h)I_{0,2},
\end{aligned}
\end{eqnarray*}
where $\alpha^{(k-s)}(h)$, $\beta^{(k-s)}(h)$, $\gamma^{(k-s)}(h)$ and $\delta^{(k-s)}(h)$ (s=1,2) are polynomials in $h$ satisfy
$$\deg\alpha^{(k-s)}(h),\ \deg\beta^{(k-s)}(h)\leq k-s-1,\ \deg\gamma^{(k-s)}(h)\leq k-s-3,\ s=1,2.$$
Noting that $I_{0,2}$ only appears in $I_{-1,2m}$, $I_{0,2m}$ and $I_{1,2l}$ and the coefficients of $I_{0,2}$ in $I_{-1,2m}$ and $I_{0,2m}$ are at most quadratic and linear polynomials of $h$, respectively, the coefficients of $I_{0,2}$ in $I_{1,2l}$ are constants, $m\geq2$, $l\geq1$. Therefore,
$$\deg\alpha^{(k)}(h),\ \deg\beta^{(k)}(h)\leq k-1,\ \deg\gamma^{(k)}(h)\leq k-3,\ \deg\delta^{(k)}(h)\leq2.$$

(ii). For $S_2$, from (2.10) and (2.11) we obtain
\begin{eqnarray}
I_{i,j}=\frac{2j}{i+j-2}(I_{i+1,j-2}-I_{i,j-2})
\end{eqnarray}
and
\begin{eqnarray}
I_{i,j}=\frac{1}{2(h-1)}(I_{i-2,j+2}+2I_{i-2,j}-4I_{i-1,j}).
\end{eqnarray}
Taking $(i,j)=(0,2),(-1,3)$ in (2.10), we get
\begin{eqnarray}
I_{0,0}=I_{1,0},\ \ I_{0,1}=I_{-1,1}.
\end{eqnarray}
From (2.20) we obtain
\begin{eqnarray}
I_{0,2}=2(h-1)I_{2,0}+2I_{1,0},\ I_{-1,3}=2(h-1)I_{1,1}+2I_{0,1}.
\end{eqnarray}
Taking $(i,j)=(-1,2)$ in (2.19), $(i,j)=(1,0)$ in (2.20) and eliminating $I_{-1,0}$, we get
\begin{eqnarray}
I_{-1,2}=\frac{4}{3}hI_{1,0}.
\end{eqnarray}
Hence, from (2.19) and (2.20) and noting that (2.21)-(2.23) we have
\begin{eqnarray}
\begin{aligned}
&I_{-1,4}=-\frac{32}{3}hI_{1,0}+8I_{0,2},\ I_{0,3}=6I_{1,1}-6I_{0,1},\ I_{3,0}=-\frac{1}{h-1}I_{1,0},\\
&I_{1,2}=\frac{2}{h-1}(I_{0,2}-2hI_{1,0}),\ I_{2,1}=\frac{1}{h-1}(I_{1,1}-2I_{0,1}),\\
&I_{-1,5}=10(4-h)I_{1,1}-40I_{0,1},\ I_{0,4}=\frac{4}{h-1}[4hI_{1,0}-(h-3)I_{0,2}],\\
&I_{1,3}=\frac{3}{h-1}[(2-h)I_{1,1}-2I_{0,1}],\ I_{2,2}=-\frac{1}{h-1}I_{0,2},\\
&I_{3,1}=\frac{1}{2(h-1)^2}(2I_{0,1}-hI_{1,1}),\ I_{4,0}=\frac{1}{(h-1)^2}I_{1,0},
\end{aligned}
\end{eqnarray}
which imply the statement (2.3) for $n=3,4$. Similar to (2.2), we can prove (2.3) by mathematical induction.

(iii). Multiplying (1.12) by $x^{i-1}y^{j-2}dx$ and integrating over $\Gamma^+_h$ yields
\begin{eqnarray}
I_{i,j}=2hI_{i+\frac{3}{2},j-2}-\frac{1}{64}I_{i+2,j-2}-\frac{1}{64}I_{i+1,j-2}.
\end{eqnarray}
If $i\neq0$ and $j\neq2$, then using (2.10) step by step, we have
\begin{eqnarray}
I_{i,j}=\begin{cases}
\sum\limits_{k=0}^{\frac{j-1}{2}}c^1_{(i,j),k}I_{i+\frac{j-1}{2}+k,1},\quad \textup{if}\ j \ \textup{odd},\\
\sum\limits_{k=0}^{\frac{j}{2}}c^2_{(i,j),k}I_{i+\frac{j}{2}+k,0},\quad\ \ \, \textup{if}\ j \ \textup{even},\\
\end{cases}
\end{eqnarray}
where $c^1_{(i,j),k}$ and $c^2_{(i,j),k}$ are arbitrary constants.
Hence, noting that $-1\leq i+\frac{j-1}{2}+k\leq n-1$ and $0\leq i+\frac{j}{2}+k\leq n$ we obtain
\begin{eqnarray}
\begin{aligned}
M(h)=&\sum\limits_{i+j=0,i\geq-1,j\geq0}^n\tilde{\rho}_{i,j}I_{i,j}(h)\\
=&\sum_{\substack{i+j=0,i\geq-1,\\j=1\,mod\,2}}^n\tilde{\rho}_{i,j}
\sum\limits_{k=0}^{\frac{j-1}{2}}c^1_{(i,j),k}I_{i+\frac{j-1}{2}+k,1}+
\sum_{\substack{i+j=0,i\geq0,\\j=0\,mod\,2,j\neq2}}^n\tilde{\rho}_{i,j}
\sum\limits_{k=0}^{\frac{j}{2}}c^2_{(i,j),k}I_{i+\frac{j}{2}+k,0}\\
&+\rho_0I_{0,2}\\
=&\sum\limits_{k=0}^n\rho^1_kI_{k-1,1}+\sum\limits_{k=0}^n\rho^2_kI_{k,0}
+\rho_0I_{0,2}\\
:=&M_1(h)+M_2(h)+\rho_0I_{0,2},
\end{aligned}
\end{eqnarray}
where $\rho^1_k$ an $\rho^2_k$ are arbitrary constants, $\rho_0=b^+_{0,2}-b^-_{0,2}-\frac{3}{4}(a^+_{1,1}-a^-_{1,1}).$

Taking $j=3$ in (2.10) and (2.25) we have
\begin{eqnarray}
I_{i,1}=\frac{1}{2i-1}\Big[2^6(4i-5)hI_{i-\frac{1}{2},1}-2(i-2)I_{i-1,1}\Big].
\end{eqnarray}
For $i\geq2$, using (2.28) step by step, we obtain
\begin{eqnarray}
I_{i,1}(h)=\alpha_{i,1}(h)I_{\frac{1}{2},1}(h)+\beta_{i,1}(h)I_{1,1}(h),
\end{eqnarray}
where $\alpha_{i,1}(h)$ and $\beta_{i,1}(h)$ are polynomials of $h$ with
$$\deg\alpha_{i,1}(h)\leq2i-3,\ \ \deg\beta_{i,1}(h)\leq2i-2.$$
For $i=0,-1$, we rewrite (2.28) as
\begin{eqnarray}
I_{i,1}=\frac{32(4i-1)}{i-1}hI_{i+\frac{1}{2},1}-\frac{2i+1}{2(i-1)}I_{i+1,1}.
\end{eqnarray}
From (2.28) and (2.30) we have
\begin{eqnarray}
\begin{aligned}
I_{0,1}=\frac{1}{2}I_{1,1}+32hI_{\frac{1}{2},1},\ I_{-\frac{1}{2},1}=64hI_{0,1},\
I_{-1,1}=80hI_{-\frac{1}{2},1}-\frac{1}{4}I_{0,1},
\end{aligned}
\end{eqnarray}
which imply
$$I_{-1,1}=(163840h^3-8h)I_{\frac{1}{2}, 1}+(2560h^2-\frac{1}{8})I_{1, 1}.$$
Hence, for $i=-1,0,1$, we have
\begin{eqnarray}
I_{i,1}=\alpha_{i,2}(h)I_{\frac{1}{2},1}(h)+\beta_{i,2}(h)I_{1,1}(h),
\end{eqnarray}
where $\alpha_{i,2}(h)$ and $\beta_{i,2}(h)$ are polynomials of $h$ with
$$\deg\alpha_{i,2}(h)\leq3,\ \ \deg\beta_{i,2}(h)\leq2.$$

Substituting (2.29) and (2.32) into $M_1(h)$ we have for $n\geq2$
\begin{eqnarray}
M_1(h)=\alpha(h)I_{\frac{1}{2},1}+\beta(h)I_{1,1},
\end{eqnarray}
where $\alpha(h)$ and $\beta(h)$ are polynomials of $h$ with $\deg\alpha(h)\leq2n-5$, $\deg\beta(h)\leq2n-4$ for $n\geq4$; $\deg\alpha(h)\leq3$, $\deg\beta(h)\leq2$ for $n=0,1,2,3$.


Taking $j=2$ in (2.10) and (2.25) we have
\begin{eqnarray}
I_{i,0}=\frac{1}{2i-3}\Big[2^8(i-2)hI_{i-\frac{1}{2},0}-(2i-5)I_{i-1,0}\Big].
\end{eqnarray}
For $i\geq1$, using (2.34) step by step, we know that $I_{i,0}$ can be written as a linear combination of $I_{0,0}$ and $I_{\frac{1}{2},0}$ with polynomial coefficients of $h$. Rewriting (2.34) as
\begin{eqnarray}
(2i-3)I_{i,0}=2^8(i-2)hI_{i-\frac{1}{2},0}-(2i-5)I_{i-1,0}
\end{eqnarray}
and taking $i=\frac{3}{2}$ in (2.35) and $i=1$ in (2.34)
we have
\begin{eqnarray}
I_{\frac{1}{2},0}=64hI_{1,0},\ \ I_{0,0}=\frac{1}{3}\big(256hI_{\frac{1}{2},0}-I_{1,0}\big).
\end{eqnarray}
Hence, we obtain for $i\geq2$
\begin{eqnarray}
I_{i,0}(h)=\gamma_{i,0}(h)I_{1,0}(h),
\end{eqnarray}
where $\gamma_{i,0}(h)$ is a polynomial of $h$ with $\deg\gamma_{i,0}(h)\leq2i-4$ and
\begin{eqnarray}
I_{0,0}(h)=\frac{1}{3}(16384h^2-1)I_{1,0}.
\end{eqnarray}
Substituting (2.37) and (2.38) into $M_2(h)$ we have for $n\geq3$
\begin{eqnarray}
M_2(h)=\gamma(h)I_{1,0}(h),
\end{eqnarray}
where $\gamma(h)$ is a polynomial of $h$ with $\deg\gamma(h)\leq2n-4$ for $n\geq3$; $\deg\gamma(h)\leq2$ for $n=0,1,2$.

 From (2.25) and (2.34) we have
$$I_{0,2}=2hI_{\frac{3}{2},0}-\frac{1}{64}I_{2,0}-\frac{1}{64}I_{1,0},\ \ I_{2,0}=I_{1,0},$$
which imply
$$I_{0,2}=2hI_{\frac{3}{2},0}-\frac{1}{32}I_{1,0}.$$
Noting that $$I_{\frac{3}{2},0}=\int_{\Gamma_h^+}xdx=\ln x_1-\ln x_2,$$
where $x_1$ and $x_2$ are the intersection points of $\Gamma_h^+$ with $x$-axis and $x_1<x_2$. By a straightforward calculation we get
$$x_1=(64h-\sqrt{4096h^2-1})^2,\ \ x_2=(64h+\sqrt{4096h^2-1})^2.$$
Hence,
\begin{eqnarray*}I_{\frac{3}{2},0}=4\ln (64h-\sqrt{4096h^2-1}),\end{eqnarray*}
which implies
\begin{eqnarray}I_{0,2}=8h\ln (64h-\sqrt{4096h^2-1})-\frac{1}{32}I_{1,0}.\end{eqnarray}

Therefore, substituting (2.33), (2.39) and (2.40) into (2.27) we obtain
\begin{eqnarray}
M(h)=\tau_3h\ln (64h-\sqrt{4096h^2-1})+\alpha_3(h)I_{\frac{1}{2},1}+\beta_3(h)I_{1,1}
+\gamma_3(h)I_{1,0},
\end{eqnarray}
where $\tau_3$ is a constant, and $\alpha_3(h)$, $\beta_3(h)$ and $\gamma_3(h)$ are polynomials of $h$ with
$$\begin{aligned}&\deg\alpha_3(h)\leq2n-5,\ \deg\beta_3(h),\deg\gamma_3(h)\leq2n-4,\ n\geq4,\\
&\deg\alpha_3(h)\leq3,\ \deg\beta_3(h),\deg\gamma_3(h)\leq2,\ n=0,1,2,3.\end{aligned}$$

(iv). Following the line of the proof of (iii), we can prove conclusion (iv).
This ends the proof.\quad $\lozenge$
 \vskip 0.2 true cm

\noindent
{\bf Lemma 2.2.}\, (i) {\it For $S_1$, the vector function $(I_{0,0},I_{1,1},I_{-1,1},I_{0,2})^T$ satisfies the following Picard-Fuchs equation
\begin{eqnarray}
\left(\begin{matrix}
                I_{0,0}\\
                 I_{1,1}\\
                 I_{-1,1}\\
                  I_{0,2}
                \end{matrix}\right)
=\frac{1}{2h+1}\left(\begin{matrix}
                2h(h+1)&0&0&0\\
                 0&\frac{1}{2}(2h+1)^2&-\frac{1}{2}&0\\
                 0&0&h^2+h&0\\
                  -2h(h+1)&0&0&\frac{1}{2}(2h+1)^2
                \end{matrix}\right)
                \left(\begin{matrix}
                I'_{0,0}\\
                 I'_{1,1}\\
                 I'_{-1,1}\\
                  I'_{0,2}
                \end{matrix}\right).
\end{eqnarray}}
\noindent
(ii) {\it For $S_2$, the vector function $(I_{0,1},I_{1,0},I_{1,1},I_{0,2})^T$ satisfies the following Picard-Fuchs equation
\begin{eqnarray}
\left(\begin{matrix}
                I_{0,1}\\
                 I_{1,0}\\
                 I_{1,1}\\
                  I_{0,2}
                \end{matrix}\right)
=\left(\begin{matrix}
                h&0&0&0\\
                 0&2h&0&0\\
                 2&0&2(h-1)&0\\
                  0&4h&0&h-1
                \end{matrix}\right)
                \left(\begin{matrix}
                I'_{0,1}\\
                 I'_{1,0}\\
                 I'_{1,1}\\
                  I'_{0,2}
                \end{matrix}\right).
\end{eqnarray}}
\noindent
(iii) {\it For $(r19)$, the vector function $(I_{\frac{1}{2},1},I_{1,1},I_{1,0})^T$ satisfies the following Picard-Fuchs equation
\begin{eqnarray}
\left(\begin{matrix}
                I_{\frac{1}{2},1}\\
                 I_{1,1}\\
                 I_{1,0}
                \end{matrix}\right)
=\left(\begin{matrix}
                h&-\frac{1}{2^6}&0\\
                 -\frac{1}{2^6}&h&0\\
                 0&0&h-\frac{1}{4096h}
                \end{matrix}\right)
                \left(\begin{matrix}
                I'_{\frac{1}{2},1}\\
                 I'_{1,1}\\
                 I'_{1,0}
                \end{matrix}\right).
\end{eqnarray}}
\noindent
(ii) {\it For $(r20)$, the vector function $(I_{0,1},I_{\frac{1}{2},1},I_{0,0})^T$ satisfies the following Picard-Fuchs equation
\begin{eqnarray}
\left(\begin{matrix}
                I_{0,1}\\
                 I_{\frac{1}{2},1}\\
                 I_{0,0}
                \end{matrix}\right)
=\left(\begin{matrix}
                h&-\frac{1}{2^6}&0\\
                 -\frac{1}{2^6}&h&0\\
                 0&0&h-\frac{1}{4096h}
                \end{matrix}\right)
                \left(\begin{matrix}
                I'_{0,1}\\
                 I'_{\frac{1}{2},1}\\
                 I'_{0,0}
                \end{matrix}\right).
\end{eqnarray}}
 \vskip 0.2 true cm

\noindent
{\bf Proof.}\, From (1.12) we get
$$\frac{\partial y}{\partial h}=\frac{x^k}{y},$$
which implies
\begin{eqnarray}
I'_{i,j}=j\int_{\Gamma^+_h}\frac{x^{i-1}y^{j-1}}{y}dx.
\end{eqnarray}
Hence,
\begin{eqnarray}
I_{i,j}=\frac{1}{j+2}I'_{i-k,j+2},\ i\geq-1,\ j\geq0.
\end{eqnarray}
Multiplying both side of (2.46) by $h$, we have
\begin{eqnarray}
hI'_{i,j}=\frac{j}{2(j+2)}I'_{i-k,j+2}+\lambda_2I'_{i-k+2,j}+\lambda_1I'_{i-k+1,j}+\lambda_0I'_{i-k,j}.
\end{eqnarray}

From (2.6) and (2.46) we have for $i\geq-1$ and $j\geq1$
\begin{eqnarray}
\begin{aligned}
I_{i,j}=&\int_{\Gamma^+_h}x^{i-k-1}y^jdx=-\frac{j}{i-k}\int_{\Gamma^+_h}x^{i-k}y^{j-1}dy\\
=&-\frac{j}{i-k}\int_{\Gamma^+_h}x^{i-k}y^{j-1}\frac{khx^{k-1}-2\lambda_2x-\lambda_1}{y}dx\\
=&-\frac{k}{i-k}hI'_{i,j}+\frac{2\lambda_2}{i-k}I'_{i-k+2,j}+\frac{\lambda_1}{i-k}I'_{i-k+1,j}.
\end{aligned}
\end{eqnarray}
(2.47)-(2.49) imply for $i\geq-1$, $j\geq1$
\begin{eqnarray}
I_{i,j}=\frac{1}{i+j-k}\big[(2-k)hI'_{i,j}-\lambda_1I'_{i-k+1,j}-2\lambda_0I'_{i-k,j}\big]
\end{eqnarray}
and
\begin{eqnarray}
I_{i,j}=\frac{2}{2i+j-2k}\big[(1-k)hI'_{i,j}+\lambda_1I'_{i-k+2,j}-\lambda_0I'_{i-k,j}\big]
\end{eqnarray}

(i) For $S_1$, from (2.47) we obtain
\begin{eqnarray}
I_{0,0}=\frac{1}{2}I'_{-1,2},\ \ I_{1,0}=\frac{1}{2}I'_{0,2},
\end{eqnarray}
and noting that (2.15) and (2.16) we have
\begin{eqnarray}
I_{0,0}=\frac{2}{2h+1}(h^2+h)I'_{0,0}.
\end{eqnarray}
From (2.49) and (2.50) we have
\begin{eqnarray}
\begin{aligned}
&I_{0,1}=\frac{1}{2}(2h+1)I'_{0,1}-\frac{1}{2}I'_{1,1},\ I_{-1,1}=\frac{1}{4}(2h+1)I'_{-1,1}-\frac{1}{4}I'_{0,1},\\
&I_{1,1}=\frac{1}{2}(2h+1)I'_{1,1}-\frac{1}{2}I'_{0,1},\ I_{0,2}=\frac{1}{2}(2h+1)I'_{0,2}-\frac{1}{2}I'_{-1,2},
\end{aligned}
\end{eqnarray}
and noting that (2.15) and (2.16) we obtain the conclusion (i).
\vskip 0.2 true cm

(ii) For $S_2$, from (2.47) and (2.19) we have
\begin{eqnarray*}
I_{1,0}=\frac{1}{2}I'_{-1,2}=\frac{2}{3}I_{1,0}+\frac{2}{3}hI'_{1,0},
\end{eqnarray*}
which implies
\begin{eqnarray}
I_{1,0}=2hI'_{1,0}.
\end{eqnarray}
From (2.49) we have
\begin{eqnarray}
\begin{aligned}
&I_{0,1}=(h-1)I'_{0,1}+I'_{-1,1},\ I_{1,1}=2(h-1)I'_{1,1}+2I'_{0,1},\\
&I_{0,2}=(h-1)I'_{0,2}+I'_{-1,2}
\end{aligned}
\end{eqnarray}
and noting that (2.21)-(2.23) we obtain the conclusion (ii).

(iii) From (2.50) and (2.51) we have
\begin{eqnarray}
I_{1,1}=hI'_{1,1}-\frac{1}{2^6}I'_{\frac{1}{2},1},\ \
I_{\frac{1}{2},1}=hI'_{\frac{1}{2},1}-\frac{1}{2^6}I'_{1,1}.
\end{eqnarray}
From (2.47) we have
\begin{eqnarray}
I_{\frac{1}{2},0}=\frac{1}{2}I'_{-1,2}.
\end{eqnarray}
From (2.10) and (2.25) we have
\begin{eqnarray}
I_{i,j}=I_{i-1,j}+\frac{2^6(4i+3j-8)}{j+2}I_{i-2,j+2}
\end{eqnarray}
and
\begin{eqnarray}
I_{i,j}=2hI_{i+\frac{3}{2},j-2}-\frac{1}{2^6}I_{i+2,j-2}-\frac{1}{2^6}I_{i+1,j-2},
\end{eqnarray}
respectively.
From (2.59) and (2.60) we get
\begin{eqnarray}
I_{i,j}=\frac{j}{2i+2j-3}\big(hI_{i+\frac{3}{2},j-2}-\frac{1}{2^6}I_{i+1,j-2}\big),
\end{eqnarray}
which implies
\begin{eqnarray}
I_{-1,2}=-2hI_{\frac{1}{2},0}+\frac{1}{32}I_{0,0}.
\end{eqnarray}
Differentiating (2.62) with respect to $h$ and noting that (2.58), we obtain
\begin{eqnarray}
I_{\frac{1}{2},0}=-\frac{h}{2}I'_{\frac{1}{2},0}+\frac{1}{128}I'_{0,0}.
\end{eqnarray}
From (2.36) and (2.63) we get $$I_{1,0}=\Big(h-\frac{1}{4096h}\Big)I'_{1,0}.$$

Following the line of the proof of (iii), we can prove the conclusion (iv).
This ends the proof.\quad $\lozenge$

 \vskip 0.2 true cm

\noindent
{\bf Lemma 2.3.}\, {\it Suppose that $h\in\Sigma$. }

\noindent
(i) {\it For $S_1$,
\begin{eqnarray}
\begin{aligned}
&I_{0,0}(h)=c_1\sqrt{h^2+h},\ \ I_{1,1}(h)=c_2h,\ \ I_{-1,1}(h)=c_2(h^2+h),\\ &I_{0,2}(h)=c_1\Big[\frac{1}{2}(2h+1)\ln|2\sqrt{h^2+h}+2h+1|-\sqrt{h^2+h}\Big].
\end{aligned}
\end{eqnarray}}
\noindent
(ii) {\it For $S_2$,
\begin{eqnarray}
\begin{aligned}
&I_{0,1}(h)=\tilde{c}_1h,\ \ I_{1,0}(h)=\tilde{c}_2\sqrt{h},\ \ I_{1,1}(h)=\tilde{c}_1-\tilde{c}_1\sqrt{1-h},\\ &I_{0,2}(h)=\tilde{c}_2\Big[\sqrt{h}-\frac{1}{2}(1-h)\ln\frac{1+\sqrt{h}}{1-\sqrt{h}}\Big].
\end{aligned}
\end{eqnarray}}

\noindent
(iii) {\it For $(r19)$,
\begin{eqnarray}
\begin{aligned}
I_{\frac{1}{2},1}(h)=\bar{c}_1\big(h-\frac{1}{2^6}\big),\ \ I_{1,1}(h)=\bar{c}_2\big(h-\frac{1}{2^6}\big),\ \ I_{1,0}(h)=\bar{c}_3\sqrt{4096h^2-1}.
\end{aligned}
\end{eqnarray}}

\noindent
(iv) {\it For $(r20)$,
\begin{eqnarray}
\begin{aligned}
I_{0,1}(h)=\hat{c}_1\big(h-\frac{1}{2^6}\big),\ \ I_{\frac{1}{2},1}(h)=\hat{c}_2\big(h-\frac{1}{2^6}\big),\ \ I_{0,0}(h)=\hat{c}_3\sqrt{4096h^2-1}.
\end{aligned}
\end{eqnarray}
where $\delta_i$ and $\tilde{\delta}_i$ ($i=1,2$) are nonzero constants.}
\vskip 0.2 true cm

\noindent
{\bf Proof.} We only prove the statement (i). The others can be shown in a similar way. By some straightforward calculations, we have $$I_{0,0}(h)=c_1\sqrt{h^2+h},\ \ I_{-1,1}(h)=c_2(h^2+h),$$
where $c_1$ and $c_2$ are nonzero constants. Therefore, we have for $h\in(0,+\infty)$
\begin{eqnarray}
\begin{aligned}
&I_{1,1}(h)=c_3(2h+1)-\frac{1}{2}c_2,\\ &I_{0,2}(h)=c_4(2h+1)+c_1\Big[\frac{1}{2}(2h+1)\ln\big(\sqrt{h^2+h}+h+\frac{1}{2}\big)-\sqrt{h^2+h}\Big],
\end{aligned}
\end{eqnarray}
where $c_3$ and $c_4$ are nonzero constants. Since $I_{1,1}(0)=I_{0,2}(0)=0$, we have $c_3=\frac{1}{2}c_2$ and $c_4=\frac{1}{2}c_1\ln2$. Substituting them into (2.68), we obtain for $h\in(0,+\infty)$
\begin{eqnarray*}
\begin{aligned}
I_{1,1}(h)=c_2h,\  I_{0,2}(h)=c_1\Big[\frac{1}{2}(2h+1)\ln(2\sqrt{h^2+h}+2h+1)-\sqrt{h^2+h}\Big].
\end{aligned}
\end{eqnarray*}
For $h\in(-\infty,-1)$, we can prove the conclusion in a similar way. This ends the proof.\quad $\lozenge$

\section{Proof of the Theorems 1.1 and 1.2}
 \setcounter{equation}{0}
\renewcommand\theequation{3.\arabic{equation}}

In the following we denote by $P_k(h)$ polynomials of $h$ with degree at most $k$.
\vskip 0.2 true cm

\noindent
{\bf Proof of the Theorem 1.1.} (1) For $S_1$ and $h\in(0,+\infty)$, let
\begin{eqnarray}
\begin{aligned}
M_1(h)=\alpha_1(h) I_{0,0}(h)+\beta_1(h) I_{1,1}(h)+\gamma_1(h) I_{-1,1}(h)+\delta_1(h) I_{0,2}(h),
\end{aligned}
\end{eqnarray}
then $M(h)$ and $M_1(h)$ have the same number of zeros in $(0,+\infty)$. From Lemmas 2.1-2.3, we have
\begin{eqnarray*}
\begin{aligned}
M_1(h)=&\alpha_1(h) I_{0,0}(h)+\beta_1(h) I_{1,1}(h)+\gamma_1(h) I_{-1,1}(h)+\delta_1(h) I_{0,2}(h),\\
=&c_1\alpha_1(h)\sqrt{h^2+h}+c_2\beta_1(h)h+c_2\gamma_1(h)(h^2+h)\\
&+c_1\delta_1(h)\Big[\frac{1}{2}(2h+1)\ln(2\sqrt{h^2+h}+2h+1)-\sqrt{h^2+h}\Big]\\
:=&P_n(h)+P_{n-1}(h)\sqrt{h^2+h}+P_3(h)\ln(2\sqrt{h^2+h}+2h+1).
\end{aligned}
\end{eqnarray*}
 After a detailed computation, we get
 \begin{eqnarray*}
\begin{aligned}
M_2(h)=\frac{d}{dh}\Big(\frac{M_1(h)}{P_3(h)}\Big)=\frac{P_{n+4}(h)+
P_{n+3}(h)\sqrt{h^2+h}}{P_6(h)\big[2h^2+2h+(2h+1)\sqrt{h^2+h}\big]}.
\end{aligned}
\end{eqnarray*}
Let $P_{n+4}(h)+P_{n+3}(h)\sqrt{h^2+h}=0$, that is,
\begin{eqnarray}P_{n+3}(h)\sqrt{h^2+h}=-P_{n+4}(h).\end{eqnarray}
By squaring the above equation, we can deduce that the number of zeros of $P_{n+4}(h)+P_{n+3}(h)\sqrt{h^2+h}$ in $(0,+\infty)$ is at most $2n+8$. Hence, noting that $2h^2+2h+(2h+1)\sqrt{h^2+h}$ does not vanish in $(0,+\infty)$, we obtain that $M_1(h)$ has at most $2n+15$ zeros in $(0,+\infty)$. Hence, $M(h)$ has at most $2n+15$ zeros for $h\in(0,+\infty)$.

In a similar way, we can prove that $M(h)$ has at most $2n+15$ zeros in $(-\infty,-1)$. Therefore, $M(h)$ has at most $4n+30$ zeros for $h\in(-\infty,-1)\cup(0,+\infty)$.

(2) For $S_2$ and $h\in(0,1)$, let
\begin{eqnarray}
\begin{aligned}
\widetilde{M}_1(h)=\alpha_2(h) I_{0,1}(h)+\beta_2(h) I_{1,0}(h)+\gamma_2(h) I_{1,1}(h)+\delta_2(h) I_{0,2}(h),
\end{aligned}
\end{eqnarray}
then $M(h)$ and $\widetilde{M}_1(h)$ have the same number of zeros in $(0,1)$.  From Lemmas 2.1-2.3, we have
\begin{eqnarray*}
\begin{aligned}
\widetilde{M}_1(h)=&\alpha_2(h) I_{0,1}(h)+\beta_2(h) I_{1,0}(h)+\gamma_2(h) I_{1,1}(h)+\delta_2(h) I_{0,2}(h),\\
=&\tilde{c}_1\alpha_2(h)h+\tilde{c}_2\beta_2(h)\sqrt{h}+\tilde{c}_1\gamma_2(h)(1-\sqrt{1-h})\\
&+\tilde{c}_2\delta_2(h)\Big[\sqrt{h}-\frac{1}{2}(1-h)\ln\frac{1+\sqrt{h}}{1-\sqrt{h}}\Big]\\
:=&P_{n-1}(h)+P_{n-1}(h)\sqrt{h}+P_{n-1}(h)\sqrt{1-h}+P_{n-1}(h)\ln\frac{1+\sqrt{h}}{1-\sqrt{h}}.
\end{aligned}
\end{eqnarray*}
 After a detailed computation, we get
 \begin{eqnarray*}
\begin{aligned}
\widetilde{M}_2(h)=&\frac{d}{dh}\Big(\frac{\widetilde{M}_1(h)}{P_{n-1}(h)}\Big)\\
=&\frac{P_{2n-1}(h)+P_{2n-1}(h)\sqrt{h}+P_{2n-1}(h)\sqrt{1-h}+P_{2n-2}(h)\sqrt{h-h^2}}
{P_{2n-2}(h)\sqrt{h-h^2}(1-\sqrt{h})^2}.
\end{aligned}
\end{eqnarray*}
Let $$P_{2n-1}(h)+P_{2n-1}(h)\sqrt{h}+P_{2n-1}(h)\sqrt{1-h}+P_{2n-2}(h)\sqrt{h-h^2}=0,$$ that is,
\begin{eqnarray}P_{2n-1}(h)\sqrt{1-h}+P_{2n-2}(h)\sqrt{h-h^2}=-P_{2n-1}(h)-P_{2n-1}(h)\sqrt{h}.\end{eqnarray}
By squaring the equation (3.4), we get
$$P_{4n-1}(h)=P_{4n-2}(h)\sqrt{h}.$$
By squaring the above equation, we can deduce that the number of zeros of $P_{2n-1}(h)+P_{2n-1}(h)\sqrt{h}+P_{2n-1}(h)\sqrt{1-h}+P_{2n-2}(h)\sqrt{h-h^2}$ in $(0,1)$ is at most $8n-2$. Noting that $\sqrt{h-h^2}(1-\sqrt{h})^2$ does not vanish in $(0,1)$, we obtain that $\widetilde{M}_1(h)$ has at most $10n-4$ zeros in $(0,1)$ for $n\geq1$. Hence, ${M}(h)$ has at most $10n-4$ zeros in $(0,1)$ for $n\geq1$. It is easy to get that ${M}(h)$ has at most 2 zeros in $(0,1)$ for $n=0$.

(3) For $(r19)$, $h\in(\frac{1}{2^6},+\infty)$ and $n\geq4$, from Lemmas 2.1 and 2.3 we obtain
\begin{eqnarray}
\begin{aligned}
M(h)=&\tau_1h\ln (64h-\sqrt{4096h^2-1})+\alpha_1(h) I_{\frac{1}{2},1}(h)+\beta_1(h) I_{1,1}(h)+\gamma_1(h)I_{1,0}(h)\\
=&\tau_1h\ln (64h-\sqrt{4096h^2-1})+\alpha_1(h)\delta_1\big(h-\frac{1}{2^6}\big)\\&+
\beta_1(h)\delta_2\big(h-\frac{1}{2^6}\big)+\gamma_1(h)\delta_3\sqrt{4096h^2-1}\\
:=&\tau_1h\ln (64h-\sqrt{4096h^2-1})+P_{2n-3}(h)+P_{2n-4}(h)\sqrt{4096h^2-1}.
\end{aligned}
\end{eqnarray}
After a detailed computation, we get
 \begin{eqnarray*}
\begin{aligned}
\widehat{M}(h)=\frac{d}{dh}\Big(\frac{M(h)}{\tau_1h}\Big)
=\frac{P_{2n-2}(h)+P_{2n-3}(h)\sqrt{4096h^2-1}}{h^2\sqrt{4096h^2-1}}.
\end{aligned}
\end{eqnarray*}
Let $P_{2n-2}(h)+P_{2n-3}(h)\sqrt{4096h^2-1}=0$, that is,
\begin{eqnarray}
P_{2n-3}(h)\sqrt{4096h^2-1}=-P_{2n-2}(h).
\end{eqnarray}
By squaring the above equation, we can deduce that the number of zeros of $P_{2n-2}(h)+P_{2n-3}(h)\sqrt{4096h^2-1}$ in $(\frac{1}{2^6},+\infty)$ is at most $4n-4$. Noting that $h^2\sqrt{4096h^2-1}$ does not vanish in $(\frac{1}{2^6},+\infty)$, we obtain that $M(h)$ has at most $4n-3$ zeros in $(\frac{1}{2^6},+\infty)$. Following the line of the above discussion, we obtain that $M(h)$ has at most 11 zeros in $(\frac{1}{2^6},+\infty)$ for $n=0,1,2,3$.

(4) For $(r20)$ and $h\in(\frac{1}{2^6},+\infty)$, similar to $(r19)$, we can prove that the number of zeros of $M(h)$ is at most $4n+3$ (reps. 8) in $(\frac{1}{2^6},+\infty)$ for $n\geq3$ (resp. $n=0,1,2$).
This end the proof of Theorem 1.1. \quad $\lozenge$
\vskip 0.2 true cm

\noindent
{\bf Proof of the Theorem 1.2.}
If $a^+_{i,j}=a^-_{i,j}$ and $b^+_{i,j}=b^-_{i,j}$, that is, the systems (1.13)-(1.16) are smooth. Since $\Gamma_h$ is symmetric with respect to $x$-axis for $h\in\Sigma$, $A_{i,2l}(h)=\oint_{\Gamma_h}x^{i-k-1}y^{2l}dx=0$, $l=0,1,2,\cdots$, where
$$\Gamma_h=\Gamma^+_h\cup\Gamma^-_h,\ \ A_{i,j}(h)=I_{i,j}(h)+J_{i,j}(h).$$
For $S_1$, $h\in(0,+\infty)$ or $h\in(-\infty,-1)$,
\begin{eqnarray*}
\begin{aligned}
M(h)=\frac{1}{2h+1}[\beta_1(h) I_{1,1}(h)+\gamma_1(h) I_{-1,1}(h)]
:=P_{n}(h),
\end{aligned}
\end{eqnarray*}
which implies that $M(h)$ has at most $n$ zeros for $h\in(0,+\infty)$  or $h\in(-\infty,-1)$. Hence, $M(h)$ has at most $2n$ zeros for $h\in(-\infty,-1)\cup(0,+\infty)$.

In a similar way, we can prove that $M(h)$ has at most $2n-1$ (resp. 1) in (0,1) if $n\geq1$ (resp. $n=0$) for $S_2$,
$M(h)$ has at most $2n-3$ (resp. 4) in $(\frac{1}{2^6},+\infty)$ if $n\geq4$ (resp. $n=0,1,2,3$) for $(r19)$ and $M(h)$ has at most $2n$ (resp. 3) in $(\frac{1}{2^6},+\infty)$ if $n\geq3$ (resp. $n=0,1,2$) for $(r20)$.

\section{Proof of the Theorem 1.3}
 \setcounter{equation}{0}
\renewcommand\theequation{4.\arabic{equation}}

From (2.8) and noting that $n=2$, we have
 \begin{eqnarray}
\begin{aligned}
M(h)=&\sum\limits_{i+j=0}^2b^+_{i,j}\int_{\Gamma^+_h}x^{i-k-1}y^jdx+
\sum\limits_{i+j=0}^2\frac{i-k-1}{j+1}a^+_{i,j}\int_{\Gamma^+_h}x^{i-k-2}y^{j+1}dx\\
&+\sum\limits_{i+j=0}^2b^-_{i,j}\int_{\Gamma^-_h}x^{i-k-1}y^jdx+
\sum\limits_{i+j=0}^2\frac{i-k-1}{j+1}a^-_{i,j}\int_{\Gamma^-_h}x^{i-k-2}y^{j+1}dx\\
=&(b^+_{0,0}-b^-_{0,0})I_{0,0}+(b_{1,0}^+-b_{1,0}^-)I_{1,0}+\big[b_{0,1}^++b^-_{0,1}+k(a_{1,0}^++a^-_{1,0})\big]I_{0,1}\\
&+(b_{2,0}^+-b^-_{2,0})I_{2,0}+\big[b^+_{1,1}+b^-_{1,1}+(1-k)(a^+_{2,0}+a^-_{2,0})\big]I_{1,1}\\
&+\big[b_{0,2}^+-b_{0,2}^--\frac{k}{2}(a_{1,1}^+-a_{1,1}^-)\big]I_{0,2}
+(k+1)(a_{0,0}^++a_{0,0}^-)I_{-1,1}\\
&-\frac{k+1}{2}(a_{0,1}^+-a_{0,1}^-)I_{-1,2}+\frac{k+1}{3}(a_{0,2}^++a_{0,2}^-)I_{-1,3}.
\end{aligned}
\end{eqnarray}

\vskip 0.2 true cm

\noindent
{\bf Lemma 4.1.}\, (i) {\it If $k=2$, $a=1$, $b=-2$ and $c=1$, then
\begin{eqnarray}
\begin{aligned}
M(h)=\sum\limits_{i=0}^5k_if_i,\ \ h\in(-\infty,-1)\cup(0,+\infty),
\end{aligned}
\end{eqnarray}
where
\begin{eqnarray*}
\begin{aligned}
&f_0=h,\ \
f_1=\sqrt{h^2+h},\ \
f_2=h^2+h,\\
&f_3=\ln|2\sqrt{h^2+h}+2h+1|,\\
&f_4=\frac{1}{2}\ln|2\sqrt{h^2+h}+2h+1|-(2h+1)\sqrt{h^2+h},\\
&f_5=\frac{1}{2}(2h+1)\ln|2\sqrt{h^2+h}+2h+1|-\sqrt{h^2+h}
\end{aligned}
\end{eqnarray*}
and
\begin{eqnarray*}
\begin{aligned}
&k_0=c_2(a_{1,0}^++a_{1,0}^--a_{0,2}^+-a_{0,2}^-+b_{1,1}^++b_{1,1}^-+b_{0,1}^++b_{0,1}^-),\\ &k_1=c_1(b_{0,0}^+-b_{0,0}^-+b_{2,0}^+-b_{2,0}^-),\\
&k_2=c_2\big[2(a_{0,0}^++a_{0,0}^-)+a_{0,2}^++a_{0,2}^-\big],\\
&k_3=\frac{c_1}{2}(b_{1,0}^+-b_{1,0}^-),\ \
k_4=\frac{c_1}{2}(a_{0,1}^+-a_{0,1}^-),\\ &k_5=c_1\big[b_{0,2}^+-b_{0,2}^--\frac{1}{2}(a_{1,1}^+-a_{1,1}^-)\big].
\end{aligned}
\end{eqnarray*}}

(ii) {\it If $k=1$, $a=\frac{1}{4}$, $b=-\frac{1}{2}$ and $c=\frac{1}{4}$, then
\begin{eqnarray}
\begin{aligned}
M(h)=\sum\limits_{i=0}^6\tilde{k}_i\tilde{f}_i,\ \ h\in(0,1),
\end{aligned}
\end{eqnarray}
where
\begin{eqnarray*}
\begin{aligned}
&\tilde{f}_0=\sqrt{h},\  \ \tilde{f}_1=h,\ \ \tilde{f}_2=h^\frac{3}{2},\\
&\tilde{f}_3=1-\sqrt{1-h},\ \
\tilde{f}_4=2h-1+(1-h)^\frac{3}{2}\\
&\tilde{f}_5=\sqrt{h}-\frac{1}{2}(1-h)\ln\frac{1+\sqrt{h}}{1-\sqrt{h}},\  \
\tilde{f}_6=\frac{\sqrt{h}-\frac{1}{2}(1-h)\ln\frac{1+\sqrt{h}}{1-\sqrt{h}}}{1-h}
\end{aligned}
\end{eqnarray*}
and
\begin{eqnarray*}
\begin{aligned}
&\tilde{k}_0=\tilde{c}_2(b_{0,0}^+-b_{0,0}^-+b_{1,0}^+-b_{1,0}^-),\quad \ \tilde{k}_2=-2\tilde{c}_2(a_{0,1}^+-a_{0,1}^-),\\ &\tilde{k}_1=\tilde{c}_1\big[b_{0,1}^++b_{0,1}^-+2(a_{1,0}^++a_{1,0}^-)+3(a_{0,0}^++a_{0,0}^-)\big],\\
&\tilde{k}_3=\tilde{c}_1\big[b_{1,1}^++b_{1,1}^--a_{2,0}^+-a_{2,0}^-)\big],\ \ \tilde{k}_4=2\tilde{c}_1(a_{0,2}^++a_{0,2}^-),\\
&\tilde{k}_5=\tilde{c}_2(b_{0,2}^+-b_{0,2}^--a_{1,1}^++a_{1,1}^-),\quad\,\tilde{k}_6=\frac{\tilde{c}_2}{2}(b_{2,0}^+-b_{2,0}^-).
\end{aligned}
\end{eqnarray*}}

\noindent
{\bf Proof.} Substituting (2.64) and (2.65) into (4.1), respectively, we get (4.2) and (4.3). This completes the proof.\quad $\lozenge$
\vskip 0.2 true cm

\noindent
{\bf Remark 4.1.} It is easy to get that
$$\frac{\partial({k}_0,{k}_1,{k}_2,{k}_3,{k}_4,k_5)}{\partial(a_{1,0}^+,b_{0,0}^+,a_{0,0}^+,b_{1,0}^+,a_{0,1}^+,b_{0,2}^+)}
=\frac{1}{2}c_1^4c_2^2\neq0,$$
$$\frac{\partial(\tilde{k}_0,\tilde{k}_1,\tilde{k}_2,\tilde{k}_3,\tilde{k}_4,\tilde{k}_5,\tilde{k}_6)}
{\partial(b_{0,0}^+,b_{0,1}^+,a_{0,1}^+,b_{1,1}^+,a^+_{0,2},b_{0,2}^+,b_{2,0}^+)}
=-2\tilde{c}_1^3\tilde{c}_2^4\neq0,$$
hence, $k_0,k_1,\cdots,k_5$ are independent and $\tilde{k}_0,\tilde{k}_1,\cdots,\tilde{k}_6$ are independent.

\vskip 0.2 true cm

In the following, we will use the {\it Chebyshev criterion} to study the number of zeros of $M(h)$ obtained in (4.2) (resp. (4.3)) for $h\in(-\infty,-1)\cup(0,+\infty)$ (resp. $h\in(0,1)$) taking into account the multiplicity. We start with some definitions and knows results.
  \vskip 0.2 true cm

\noindent
{\bf Definition 4.1.}\, \cite{MV} Let $\varphi_0,\varphi_1,\cdots,\varphi_{n-1}$ be analytic functions on an open interval $L$. $(\varphi_0,\varphi_1,\cdots,\varphi_{n-1})$ is an {\it extended complete Chebyshev system} (in short, ECT-system) on $L$ if, for all $k=1,2,\cdots,n$, any nontrivial linear combination
\begin{eqnarray}\lambda_0\varphi_0(x)+\lambda_1\varphi_1(x)+\cdots+\lambda_{k-1}\varphi_{k-1}(x)\end{eqnarray}
has at most $k-1$ isolated zeros on $L$ counting multiplicities.
  \vskip 0.2 true cm

\noindent
{\bf Remark 4.2.}\, If $(\varphi_0,\varphi_1,\cdots,\varphi_{n-1})$ is an ECT-system on $L$, then for each $k=1,2,\cdots,n$, there exists a linear combination (4.4) with $k-1$ simple zeros on $L$ (see \cite{GV})
  \vskip 0.2 true cm
In order to show that a set of functions is a ECT-system, the notion of the {\it Wronskian} proves to be extremely useful.
\vskip 0.2 true cm

\noindent
{\bf Definition 4.2.}\, \cite{MV} Let $\varphi_0,\varphi_1,\cdots,\varphi_{n-1}$ be analytic functions on an open interval $L$. The {\it Wronskian} of $(\varphi_0,\varphi_1,\cdots,\varphi_{k-1})$ at $x\in L$ is
$$W[\varphi_0,\varphi_1,\cdots,\varphi_{k-1}](x)=det\Big(\varphi_j^{(i)}(x)\Big)_{0\leq i,j\leq k-1}$$
\vskip 0.2 true cm

For the sake of shortness we will use the notation
$$W[\varphi_0,\varphi_2,\cdots,\varphi_{k-1}](x)=W_k(x).$$
The following lemma is well-known (see \cite{KS}).
 \vskip 0.2 true cm

\noindent
{\bf Lemma 4.2.}\, {\it $(\varphi_0,\varphi_1,\cdots,\varphi_{n-1})$ is an ECT-system on $L$ if and only if, for each $k=1,2,\cdots,n$,
 $$W_k(x)\neq0 \ \textup{for\ all}\ x\in L.$$}
 \vskip 0.2 true cm

 First of all, we will prove that $(f_0,f_1,f_2,f_3,f_4,f_5)$ in (4.2) is an ECT-system for $h\in(0,+\infty)$. By some straightforward calculations we have
\begin{eqnarray}
\begin{aligned}
&W_1(h)=h,\ \ W_2(h)=-\frac{h}{2\sqrt{h^2+h}},\ \ W_3(h)=-\frac{h(4h+3)}{4(1+h)\sqrt{h^2+h}},\\
&W_4(h)=-\frac{1}{4(h^2+h)^3}\big[8h^3+15h^2+6h-3(2h+1)\sqrt{h^2+h}\ln(2\sqrt{h^2+h}+2h+1)\big],\\
&W_5(h)=\frac{3}{16\sqrt[11]{(h^2+h)^2}}\big[4h^3-2h^2-6h+3\sqrt{h^2+h}\ln(2\sqrt{h^2+h}+2h+1)\big],\\
&W_6(h)=-\frac{3}{16(h^2+h)^8}\big[8h^4-8h^3-46h^2-30h\\
&\qquad\qquad +3(6h+5)\sqrt{h^2+h}\ln(2\sqrt{h^2+h}+2h+1)\big].
\end{aligned}
\end{eqnarray}
It is easy to check that $W_k(h)$ $(k=1,2,\cdots,6)$ do not vanish for $h\in(0,+\infty)$ by Matlab or Maple. The proof for $h\in(-\infty,-1)$ follows in the same way.

Secondly we begin to prove that $(\tilde{f}_0,\tilde{f}_1,\tilde{f}_2,\tilde{f}_3,\tilde{f}_4,\tilde{f}_5,\tilde{f}_6)$ in (4.3) is an ECT-system for $h\in(0,1)$. By some straightforward calculations we have
\begin{eqnarray}
\begin{aligned}
&W_1(h)=\sqrt{h},\ \ W_2(h)=\frac{1}{2}\sqrt{h},\ \ W_3(h)=\frac{1}{4},\\
&W_4(h)=\frac{3}{32h^3\sqrt[5]{(1-h)^2}}\big[4h^2-5h+2-2(1-h)^\frac{5}{2}\big],\\
&W_5(h)=\frac{9}{512h^5\sqrt[11]{(1-h)^2}}\big[5h^4-37h^3+71h^2-51h+12\\
&\qquad\qquad+(15h^3-50h^2+45h-12)(1-h)^\frac{1}{2}\big],\\
&W_6(h)=\frac{9}{8192h^\frac{15}{2}(1-h)^9}\Big[360h^4-980h^3+930h^2-330h+24\\
&\qquad\qquad+2(280h^3-375h^2+147h-12)(1-h)^\frac{3}{2}\\
&\qquad\qquad-45(4h-3)(1-h)^2h^\frac{3}{2}\ln\frac{1+\sqrt{h}}{1-\sqrt{h}}\Big],\\
&W_7(h)=\frac{27}{262144h^{13}\sqrt[33]{(1-h)^2}}\Big[16(26950h^8-26215h^7+17625h^6-49857h^5\\
&\qquad\qquad+58944h^4-31497h^3+9165h^2-1335h+60)(1-h)^2\\
&\qquad\qquad+2(137445h^9-290980h^8+217735h^7-315480h^6+620745h^5\\
&\qquad\qquad-584886h^4+285531h^3-80070h^2+11400h-480)(1-h)^\frac{1}{2}\\
&\qquad\qquad-15h^\frac{3}{2}(1-h)^\frac{5}{2}(9163h^6-4127h^5+5848h^4-16224h^3\\
&\qquad\qquad+11895h^2-4305h+630)\ln\frac{1+\sqrt{h}}{1-\sqrt{h}}\Big].\\
\end{aligned}
\end{eqnarray}
It is easy to check that $W_k(h)$ $(k=1,2,\cdots,7)$  do not vanish for $h\in(0,1)$ by Matlab or Maple.
 \vskip 0.2 true cm

\noindent
{\bf Proof of Theorem 1.3.} From Lemma 4.2 we know that $(f_0,f_1,\cdots,f_5)$ in (4.2) (resp. $(\tilde{f}_0,\tilde{f}_1,\cdots,\tilde{f}_6)$ in (4.3)) is an ECT-system for $h\in(-\infty,-1)$ or $(0,+\infty)$ (resp. $h\in(0,1)$) which implies that $M(h)$ in (4.1) has at most 5 (resp. 6) zeros in $(-\infty,-1)$ or $(0,+\infty)$ (resp. $(0,1)$), and this number is realizable by Remark 4.1. Hence, the Abelian integral $M(h)$ defined in (4.1) has exact 5 (resp. 6) zeros in $(-\infty,-1)$ or $(0,+\infty)$ (resp. $(0,1)$). This completes the proof of Theorem 1.3.\quad $\lozenge$
 \vskip 0.2 true cm

\noindent
{\bf Remark 4.3.} (i) If $a^+_{i,j}=a^-_{i,j}$ and $b^+_{i,j}=b^-_{i,j}$, then $k_1=k_3=k_4=k_5=0$ and $\tilde{k}_0=\tilde{k}_2=\tilde{k}_5=\tilde{k}_6=0$. Hence, from the above discussion we obtain that $M(h)$ has exact one zero in $(-\infty,-1)$ or $(0,+\infty)$ for $k=1$ and $M(h)$ has exact two zeros in $(0,1)$ for $k=2$ which coincide with the result of \cite{CJ}. 
 \vskip 0.2 true cm

\noindent
(ii) In the \cite{LM,LC,CLZ}, the authors study the bifurcations of limit cycle for $S_1$, $S_2$, $S_3$ and $S_4$ when they are perturbed inside discontinuous quadratic polynomials of the form (1.8) by using the averaging method of first order. But the averaged functions involve a great quantity of computations even if $n=2$. So that we can not use the method to study them for arbitrary $n$.

 \vskip 0.2 true cm

\noindent
{\bf Acknowledgment}
 \vskip 0.2 true cm

\noindent
The second author is supported by National Natural Science Foundation of China (11671040), the first author is supported by National Natural Science Foundation of China (11671040),  Science and Technology Pillar Program of Ningxia(KJ[2015]26(4)),  Visual Learning Young Researcher of Ningxia and  Key Program of Ningxia Normal University(NXSFZD1708).

\end{document}